\documentclass{article} 
\usepackage{amsfonts}

\usepackage[pdftex]{graphicx}

\usepackage{latexsym,amssymb,amsthm,amsmath}
\newtheorem{theorem}{Theorem}[section]

\theoremstyle{proposition}
\theoremstyle{definition}

\theoremstyle{remark}

\begin{document}

\title{A Natural Representation of Volumes Yields a Remarkable Affine Consequence}

\author{{Wladimir G. Boskoff,  Bogdan D. Suceav\u{a}}}
\maketitle

{\it To appear in the Houston Journal of Mathematics.}

\begin{abstract}
At the beginning of the 20th Century there was a growing interest for the investigation of the action of linear groups on the geometry of surfaces. In that context of ideas, the quest for a connection between curvature and the behaviour of linear groups rose naturally. Pursuing the original thought, we investigate how the geometric meaning of this idea is intimately related to the concept of volume of parallelepiped boxes. We show how the ratio of the Gaussian curvature divided by the fourth power of a certain distance of interest in the geometry of surfaces can be represented as a function of volumes. This geometric description explores the profound meaning of a quantity considered by {\c{T}}i{\c{t}}eica in 1907, in a work that sparked a growing interest in affine differential geometry, as an illustration of Felix Klein's Erlangen Program, in which the quest for geometric invariants was the main point of inquiry.
\end{abstract}

{\bf Keywords:} affine invariants; Erlangen program; Gaussian curvature; volume.

{\bf MSC 2020:}	53A15;  	53A05; 	53A99; 	28A75.


\section{Introduction}
In his much cited work \cite{K} published originally in 1908, Felix Klein (1849-1925) starts the investigation of the geometric representations of space by looking at the simplest manifolds. To fully support his investigation, F. Klein introduces determinants early in his exposition, and discusses the interpretation of signs (in the definition of the determinant), as well as the volume of polyhedrons. In the same historical period when Felix Klein reflected upon the first edition of his seminal volume, a former doctoral student of Gaston Darboux (1842-1917),  Gheorghe {\c{T}}i{\c{t}}eica (1873-1939), introduced \cite{T1907,T1908x,T1908b} the investigation of a new class of surfaces, which have been later called by Wilhelm Blaschke (1885-1962) and other authors {\c{T}}i{\c{t}}eica surfaces \cite{Blaschke}. An important fact is that {\c{T}}i{\c{t}}eica discovered that a whole family of surfaces remains invariant under the actions of the so-called centro-affine group. 

{\c{T}}i{\c{t}}eica's idea was an excellent illustration of Felix Klein's main point in his Erlangen Program.  In his 1872 work, Klein formulated what he described as a {\it comprehensive problem}:
``Given a manifoldness and a group of transformations of the same; to investigate the configurations belonging to the manifoldness with regard to such properties as are not altered by the
transformations of the group." This was the original quest. {\c{T}}i{\c{t}}eica's construction  \cite{T1907,T1908x,T1908b} appeared later as both an illustration and as a surprise.

While the geometric idea had much success at that time \cite{A}, what is truly interesting to see is a discussion at fundamental level of how the original affine invariant actually appeared. 
To better see this, consider an arbitrary point
$p(x,y)$ lying in the domain of a surface patch locally described in graph form as $f(x,y)=\langle x,y,u(x,y)\rangle.$
 
It is well known (see e.g. \cite{BC}) that the Gaussian curvature   $K_f(p)$ at that point has the expression
 $$K_f \left( p\right)=\dfrac{\det h_{ij}}{\det g_{ij}}=\dfrac{\dfrac{\partial^2 u}{\partial x^2}\dfrac{\partial^2 u}{\partial y^2} -\left(\dfrac{\partial^2 u}{\partial x \partial y}\right)^2 }{\left[\left( 
 	\dfrac{\partial u}{\partial x}\right) ^{2}+\left( \dfrac{\partial u}{\partial y 
 	}\right) ^{2}+1 \right]^2}.$$

\section{A fundamental approach to volumes.}
Consider the basis of the tangent plane at $f(p)$ to the surface $f$, i.e. $$\dfrac{\partial f}{\partial x}=\left<1,0,\dfrac{\partial u}{\partial x}\right>; \ \dfrac{\partial f}{\partial y}=\left<0,1,\dfrac{\partial u}{\partial y}\right>.$$

Denoting by $\vec{n}:= \dfrac{\partial f}{\partial x}\times \dfrac{\partial f}{\partial y}=\left<-\dfrac{\partial u}{\partial x}, -\dfrac{\partial u}{\partial y}, 1\right>$, we remark that the numerator in the Gaussian curvature formula represents the oriented area of the parallelogram spanned by the vectors $\dfrac{\partial \vec{n}}{\partial x}=\left<-\dfrac{\partial^2 u}{\partial x^2},-\dfrac{\partial^2 u}{\partial x \partial y},0\right>$ and  $\dfrac{\partial \vec{n}}{\partial y}=\left<-\dfrac{\partial^2 u}{\partial y\partial x},-\dfrac{\partial^2 u}{\partial y^2},0\right>,$ expressed with respect to the vector $\dfrac{\partial \vec{n}}{\partial x}\times \dfrac{\partial \vec{n}}{\partial y} $.\\

The denominator is the square of the determinant of the coefficients of the metric attached to the surface $f$, namely: $$ds^2=\left(1+\left(\dfrac{\partial u}{\partial x}\right)^2\right)dx^2 +2 \dfrac{\partial u}{\partial x}\dfrac{\partial u}{\partial y}dxdy+\left(1+\left(\dfrac{\partial u}{\partial y}\right)^2\right)dy^2.$$

To complete the full image of what is needed, consider the four vectors constructed with the second order partial derivatives, as follows. $$ \dfrac{\partial^2 f}{\partial x^2}=\left<0,0,\dfrac{\partial^2 u}{\partial x^2}\right>; \ \dfrac{\partial^2 f}{\partial y^2}=\left<0,0,\dfrac{\partial^2 u}{\partial y^2}\right>; \ \dfrac{\partial^2 f}{\partial x \partial y}=\dfrac{\partial^2 f}{\partial y \partial x}=\left<0,0,\dfrac{\partial^2 u}{\partial x \partial y}\right>.$$ 
All these vectors are parallel to the canonical basis   unit vector $\langle 0,0,1 \rangle $.\\

It would be interesting to consider three 
 parallelepipeds, specifically those determined by the vectors $\left\lbrace\dfrac{\partial^2 f}{\partial x^2},\ \dfrac{\partial f}{\partial x},\ \dfrac{\partial f}{\partial y}\right\rbrace $,  $\left\lbrace\dfrac{\partial^2 f}{\partial y^2},\ \dfrac{\partial f}{\partial x},\ \dfrac{\partial f}{\partial y}\right\rbrace$ and $\left\lbrace\dfrac{\partial ^2 f}{\partial x \partial y},\ \dfrac{\partial f}{\partial x},\ \dfrac{\partial f}{\partial y}\right\rbrace$.
Their oriented volumes are 

$$
V_x :=
\left |
\begin{array}{ccc}
0 & 0 & \dfrac{\partial^2 u}{\partial x^2} \\[2ex]
1 & 0 & \dfrac{\partial u}{\partial x } \\[2ex]
0 & 1 & \dfrac{\partial u}{\partial y} \\
\end{array}
\right |,\quad 
V_y:=\left |
\begin{array}{ccc}
0 & 0 & \dfrac{\partial^2 u}{\partial y^2} \\[2ex]
1 & 0 & \dfrac{\partial u}{\partial x } \\[2ex]
0 & 1 & \dfrac{\partial u}{\partial y} \\
\end{array}
\right |, \quad
V_{xy}:=\left |
\begin{array}{ccc}
0 & 0 & \dfrac{\partial^2 u}{\partial x \partial y} \\[2ex]
1 & 0 & \dfrac{\partial u}{\partial x } \\[2ex]
0 & 1 & \dfrac{\partial u}{\partial y} \\
\end{array}
\right |. 
$$
Remark that the numerator of the Gaussian curvature (previously expressed as an oriented area) is now described by the equality 
$$ \frac{\partial^2 u}{\partial x^2}\frac{\partial^2 u}{\partial y^2} -\left(\frac{\partial^2 u}{\partial x \partial y} \right)^2=V_xV_y-(V_{xy})^2,$$ therefore $$K_f(p)=\dfrac{V_xV_y-(V_{xy})^2}{\left[ \left(\dfrac{\partial u}{\partial x}\right)^2+\left(\dfrac{\partial u}{\partial y}\right)^2  +1\right]^2}.$$

\section{The path to {\c{T}}i{\c{t}}eica's Theorem}
The part explored in the previous section, i.e. $K_f(p)$, is the first of the two components investigated by {\c{T}}i{\c{t}}eica \cite{T1907,T1908x,T1908b}. Now we will approach the second component, which is considered as follows. 
Denote by $d_f(p)$ the distance from the origin $O(0,0,0)$ to the tangent plane to the surface at the point $f(p)$,
$$
d_{f}(p) =\dfrac{\left|x\dfrac{\partial u
	}{\partial x}+y\dfrac{\partial u}{\partial y} - u\left( x,y\right)\right|  }{\sqrt {\left( 
		\dfrac{\partial u}{\partial x}\right) ^{2}+\left( \dfrac{\partial u}{\partial y 
		}\right) ^{2}+1 }}.$$

Consider the parallelepiped generated by the vectors $\left\lbrace f,\ \dfrac{\partial f}{\partial x},\ \dfrac{\partial f}{\partial y}\right\rbrace$. Its oriented volume is expressed by the determinant 
$$
V:=\left |
\begin{array}{ccc}
x & y & u(x,y) \\[2ex]
1 & 0 & \dfrac{\partial u}{\partial x } \\[2ex]
0 & 1 & \dfrac{\partial u}{\partial y} \\
\end{array}
\right |.  
$$

The really striking coincidence of components yields the remarkable relation
  $$V^4=\left[ \left(\dfrac{\partial u}{\partial x}\right)^2+\left(\dfrac{\partial u}{\partial y}\right)^2  +1\right]^2\cdot d^4_f(p),$$ 
i.e.
$$\dfrac{K_f(p)}{d^4_f(p)}=\dfrac{V_xV_y-(V_{xy})^2}{V^4}.$$
It follows that the numerator and the denominator of the previous ratio  are respectively
$$
\left |
\begin{array}{ccc}
0 & 0 & \dfrac{\partial^2 u}{\partial x^2}\vphantom{\dfrac{\partial^2 u}{\partial y \partial x}} \\[2ex]
1 & 0 & \dfrac{\partial u}{\partial x } \\[2ex]
0 & 1 & \dfrac{\partial u}{\partial y} \\
\end{array}
\right |\cdot 
\left |
\begin{array}{ccc}
0 & 0 & \dfrac{\partial^2 u}{\partial y^2} \\[2ex]
1 & 0 & \dfrac{\partial u}{\partial x } \\[2ex]
0 & 1 & \dfrac{\partial u}{\partial y} \\
\end{array}
\right | -
\left |
\begin{array}{ccc}
0 & 0 & \dfrac{\partial^2 u}{\partial x \partial y} \\[2ex]
1 & 0 & \dfrac{\partial u}{\partial x } \\[2ex]
0 & 1 & \dfrac{\partial u}{\partial y} \\
\end{array}
\right |\cdot
\left |
\begin{array}{ccc}
0 & 0 & \dfrac{\partial^2 u}{\partial y \partial x} \\[2ex]
1 & 0 & \dfrac{\partial u}{\partial x } \\[2ex]
0 & 1 & \dfrac{\partial u}{\partial y} \\
\end{array}
\right | 
$$
and 

$$
\left |
\begin{array}{cccc}
x & y & u(x,y) \\[2ex]
1 & 0 & \dfrac{\partial u}{\partial x } \\[2ex]
0 & 1 & \dfrac{\partial u}{\partial y} \\
\end{array}
\right |^4.  
$$

 A centro-affine  transformation of the surface $f$ consists of the product between the vector representation of the surface $f$ and a $3 \times 3$ matrix $\mathbb{A}=(a_{ij})$, with $\det \mathbb{A} \neq 0$.\\ 

$$\bar{f}(x,y)=f(x,y)\cdot \mathbb{A}= $$ $$=\langle a_{11}x+a_{21}y+a_{31} u(x,y), a_{12}x+a_{22}y+a_{32}u(x,y), a_{13}x+a_{23}y+a_{33}u(x,y)\rangle$$

 Denote by $\bar{f}(q)=f(p)\cdot \mathbb{A}.$
For $\bar{f},$ the vector basis we are using is 
$$\frac{\partial \bar f}{\partial x}=\left<a_{11}+a_{31}\frac{\partial u}{\partial x}, a_{12}+a_{32}\frac{\partial u}{\partial x} , a_{13}+a_{33}\frac{\partial u}{\partial x}  \right>$$  
$$\frac{\partial \bar f}{\partial y}=\left<a_{21}+a_{31}\frac{\partial u}{\partial y}, a_{22}+a_{32}\frac{\partial u}{\partial y} , a_{23}+a_{33}\frac{\partial u}{\partial y}  \right> $$ and the four other vectors are
$$\frac{\partial^2 \bar f}{\partial x^2}=\frac{\partial^2 u}{\partial x^2}\cdot \langle a_{31}, a_{32}, a_{33}\rangle, \ \frac{\partial^2 \bar f}{\partial x \partial y}=\frac{\partial^2 \bar f}{\partial y \partial x}=\frac{\partial^2 u}{\partial x \partial y}\cdot \langle a_{31}, a_{32}, a_{33}\rangle,$$ $$ \ \frac{\partial^2 \bar f}{\partial y^2}=\frac{\partial^2 u}{\partial y^2}\cdot \langle a_{31}, a_{32}, a_{33}\rangle.$$
The first determinant of the numerator is
$$
\left |
\begin{array}{cccc}
a_{31}\dfrac{\partial^2 u}{\partial x^2} & a_{32}\dfrac{\partial^2 u}{\partial x^2} & a_{33}\dfrac{\partial^2 u}{\partial x^2} \\[2ex]
a_{11}+a_{31}\dfrac{\partial u}{\partial x } & a_{12}+a_{32}\dfrac{\partial u}{\partial x } & a_{13}+a_{33}\dfrac{\partial u}{\partial x } \\[2ex]
a_{21}+a_{31}\dfrac{\partial u}{\partial y } & a_{22}+a_{32}\dfrac{\partial u}{\partial y } & a_{23}+a_{33}\dfrac{\partial u}{\partial y} \\
\end{array}
\right |=\det \mathbb{A} \cdot \dfrac{\partial^2 u}{\partial x^2}.
$$
Therefore we see that the new numerator has the value $$\bar V_{\bar x}\bar V_{\bar y}-\bar V^2_{\bar x \bar y}=\det \mathbb{A}^2 \cdot \left(\dfrac{\partial^2 u}{\partial x^2}\dfrac{\partial^2 u}{\partial y^2} -\left(\dfrac{\partial^2 u}{\partial x \partial y} \right)^2\right). $$
The denominator is described by the $4^{th}$ power of a determinant, namely
$$
\bar V=\left |
\begin{array}{ccc}
a_{11}x+a_{21}y+a_{31}u & a_{12}x+a_{22}y+a_{32}u & a_{13}x+a_{23}y+a_{33}u \\[2ex]
a_{11}+a_{31}\dfrac{\partial u}{\partial x } & a_{12}+a_{32}\dfrac{\partial u}{\partial x } & a_{13}+a_{33}\dfrac{\partial u}{\partial x } \\[2ex]
a_{21}+a_{31}\dfrac{\partial u}{\partial y } & a_{22}+a_{32}\dfrac{\partial u}{\partial y } & a_{23}+a_{33}\dfrac{\partial u}{\partial y} \\
\end{array}
\right |.
$$
We observe that $\bar{V}=V\cdot \det \mathbb{A},$ i.e. the denominator has the value $$\det \mathbb{A}^4 \cdot \left(x\dfrac{\partial u}{\partial x}+y%
\dfrac{\partial u}{\partial y}-u\left( x,y\right)\right)^4. $$ 
Therefore, we proved the following  \cite{T1907}.
	
\begin{theorem} With the notations described previously,
$$\frac{ K_{\bar f}(q)}{\Large ({d}_{\bar f})^4(q)}=\dfrac{1}{(\det \mathbb{A})^2}\cdot\dfrac{K_f \left( p\right)}{({d_{f}})^{4}\left( p\right)}
	.    $$ 
	\end{theorem}
In terms of volumes the result can be written as
	$$\dfrac{\bar V_{\bar x}\bar V_{\bar y}-\bar V^2_{\bar x \bar y}}{\bar V^4}=\dfrac{1}{(\det \mathbb{A})^2}\cdot \dfrac{V_x V_y - V_{xy}}{V^4}$$
since we proved that the bar-volumes are related to the initial volumes through the formulas $\bar V= \det \mathbb{A} \cdot V$, etc.

Since this is a key illustration of the Erlangen Program's philosophy, it is of major interest to describe this result's importance. Basically, the theorem says that 
 the centro-affine transformations  preserve the ratio $\dfrac{K}{d^4}$ up to a constant. Thus, this invariant property responds vividly to Felix Klein's quest  for ``properties that are not altered by the transformations of the group". As a bonus, this invariance property involves the Gaussian curvature in the conversation, something that in 1907 was perceived as a novel element.

 Additionally, if we choose the matrix $\mathbb{A}$ such that the square of the determinant is $1$ the ratio $\dfrac{K}{d^4}$ is fully preserved. Following the terminology rooted in the studies written after WW I, 
the surface is called a {\c{T}}i{\c{t}}eica surface if at each point $p$ the ratio $\dfrac{K_f(p)}{d^4_f(p)}$ is a constant, let us denote it by $R_f$.
Normally $R_{\bar{f}} =\dfrac{1}{(\det \mathbb{A})^2}R_f$, but if $(\det \mathbb{A})^2=1$ we have $R_{\bar{f}}=R_f$. Therefore the very fact that a given surface is a  {\c{T}}i{\c{t}}eica surface is related to the conservation of volumes, and we illustrated above how this process takes place. \\

\section{Why geometry alone does not characterize {\c{T}}i{\c{t}}eica surfaces} 
I. The position of an object in space can affect its centro-affine structure.\\

Consider a sphere centered at origin, with radius  $R$. Its Gaussian curvature at each point is $K_f=\dfrac{1}{R^2}.$ The distance from the origin to the tangent plane is  $R$. Therefore, at each point $\dfrac{K_f}{d^4_f}=\dfrac{1}{R^6}$, i.e. the sphere centered at the origin is a {\c{T}}i{\c{t}}eica surface.\\
Consider moving the sphere somewhere else in space. Its Gaussian  curvature remains constant. But the distance from the origins to the tangent planes varies, hence such a sphere is not a {\c{T}}i{\c{t}}eica surface.

II. The property of a surface to be a {\c{T}}i{\c{t}}eica surface does not depend on the geometry of the surface.\\

To illustrate this idea, consider the pseudosphere endowed with the metric induced on it by the three-dimensional Euclidean space. Hence, its metric is $$ds^2=\sin^2 x_2                   \ dx_1^2+\cot^2 x_2\ dx_2^2.$$
As it is well-known, this metric has a constant Gaussian curvature, $K=-1.$ Furthermore, if we perform a change of coordinates 
$$
\left \{
\begin{array}{l}
x= x_1\\
y=\dfrac{1}{ \sin x_2}\ ,\\
\end{array}
\right.    
$$
then the previous metric of the pseudosphere is transformed into the Poincar\'{e} metric of the half-plane $$ds^2=\dfrac{dx^2+dy^2}{y^2}.  $$ This shows that the geometry of the pseudosphere is the classical non-Euclidean geometry. The pseudosphere is not a {\c{T}}i{\c{t}}eica surface because $\dfrac{K}{d^4} $ is not a constant. \\

The metric of the Poincar\'{e} disk $$ds^2=\dfrac{4}{(1-y_1^2-y_2^2)^2}(dy_1^2+dy_2^2)$$ and the previous Poincar\'{e} half plane are equivalent; a point $(x,y)$ in the Poincar\'{e} half-plane model maps to $$ y_1=\dfrac{2x}{x^2+(1-y)^2}  ,\  y_2=\dfrac{1-x^2-y^2}{x^2+(1-y)^2} $$ in the Poincar\'{e}  disk-model, and one metric is obtained from the other upon pull-back. Therefore they share the classical non-Euclidean geometry.\\ If $(x_0,x_1,x_2)$ is a vector of the three-dimensional coordinate space $\mathbf{R}^3$, the Minkowski quadratic form is defined by $$Q(x_0,x_1,x_2)=-x_0^2+x_1^2+x_2^2.$$ The vectors $v \in \mathbf{R}^3$ such that $Q(v)=-1$ form a two-dimensional hyperboloid consisting of the two connected components (the two "sheets"); the forward sheet, where $x_0 > 0$, and the backward sheet, where $x_0<0$. The points of the two-dimensional hyperboloid model are the points lying on the forward sheet. From the Minkowski point of vue this hyperboloid is a Minkowski sphere, its radius is $1$. The Minkowski metric attached to the Minkowski quadratic form is $$ds^2=-dx_0^2+dx_1^2+dx_2^2.$$ If $$x_0=\cosh u_1; \ x_1=\sinh u_1 \cos u_2; \ x_2=\sinh u_1 \sin u_2,$$ the metric of this Minkowski sphere is $$ds^2=du_1^2+\sinh^2 u_1 du_2^2.$$ The Gaussian curvature of this metric is $K=-1$. If we consider the change of coordinates (see \cite{BC})  $$u_1=2 \tanh^{-1}(y_1^2+y_2^2); \ u_2=\arctan\dfrac{y_2}{y_1},$$ replacing in the metric of the Minkowski sphere we obtain the metric of the Poincar\'{e} disk $$ds^2=\dfrac{4}{(1-y_1^2-y_2^2)^2}(dy_1^2+dy_2^2).$$ Therefore the hyperboloid (the Minkowski sphere), the Poincar\'{e} disk, the Poincar\'{e} half-plane  and the pseudosphere share the same classical non-Euclidian geometry. For the pseudosphere, however, the ratio  $\dfrac{K}{d^4}$ is not constant,; hence the pseudosphere is not a {\c{T}}i{\c{t}}eica surface. In contrast, for the Minkowski sphere the Minkowski distance from the origin to any tangent plane equals the radius, namely $d=1$, and consequently the ratio $\dfrac{K}{d^4}=-1$ is a constant. This shows that the Minkowski sphere is a {\c{T}}i{\c{t}}eica surface. Although the {\c{T}}i{\c{t}}eica condition is classically proven for surfaces in Euclidean space,  in \cite{boss}, the notion of {\c{T}}i{\c{t}}eica surface is extended to Minkowski spaces by direct computation. Therefore the pseudosphere and the Minkowski sphere share the same classical non-Euclidean geometry but not the {\c{T}}i{\c{t}}eica property. Let us observe that the explanation of the {\c{T}}i{\c{t}}eica invariant through oriented volumes offered in this paper can be easily extended to Minkowski surfaces replacing the direct computation from \cite{boss}.

The theme of affine differential geometry enjoyed much interest over many decades (see e.g. \cite{Calabi,Chen,Nomizu, Simon}) and a fundamental explanation of its connection with a calculus-like elementary geometric measure can only shed an additional light on such an interesting topic.

\section{Acknowledgements} The authors would like to express their gratitude to the anonymous referee for the careful reading of the manuscript and for the insightful suggestions that significantly improved the clarity and the professional appearance of this paper. In particular, we are grateful for the recommendations regarding the geometric interpretation of the Poincar\'{e} models and for the correction of the historical references.

\vspace{.2cm}

Wladimir G. Boskoff\\
Department of Mathematics and Computer Science, Ovidius University,  Bd. Mamaia, nr. 124, Constanța, Romania \\ 
email: boskoff@univ-ovidius.ro

\vspace{.3cm}

Bogdan D. Suceav\u{a}\\
Department of Mathematics, California State University, Fullerton, McCarthy Hall 154, Fullerton, CA 92831-6850 \\
email: bsuceava@fullerton.edu

\end{document}